\UseRawInputEncoding
\documentclass[10pt]{article}
 \usepackage{mathrsfs}
 \usepackage{float}
 \usepackage{amsfonts}
\pdfoutput=1
\usepackage{amssymb,amsmath,amsthm,amsfonts}
\usepackage[dvips]{graphicx}
\usepackage{subfigure}
\usepackage{color}
\usepackage{indentfirst}
\numberwithin{equation}{section}

\begin{document}

\title{\bf Stabilization of  Stackelberg Game-Based Control Systems
}

\date{}

\author{Yue Sun, Juanjuan Xu, Huanshui Zhang
\thanks{This work is supported by the National Natural Science Foundation of China under Grants 61633014, 61873332, U1806204, U1701264, 61922051, the foundation for Innovative Research Groups of National Natural Science Foundation of China (61821004). Y. Sun, J. Xu and H. Zhang are with the School of Control Science and Engineering,
Shandong University, Jinan, Shandong 250061, China (e-mail:
sunyue9603@163.com; juanjuanxu@sdu.edu.cn; hszhang@sdu.edu.cn).}%
}


\maketitle

\begin{abstract}
In this paper, we are concerned with the stabilizatbility of Stackelberg game-based systems. In particular, two players are involved in the system where one is the follower to minimize the related cost function and the other is the leader to stabilize the system. The main contribution is to derive the necessary and sufficient condition for the stabilization of the game-based system. The key technique is to explicitly solve the forward and backward difference equations (FBDEs) based on the maximum principle and give the optimal feedback gain matrix of the leader by using the matrix maximum principle.

\vspace{0.3cm}


\vspace{0.3cm}

\end{abstract}

\section{Introduction}
In modern control theory, stabilizability is an essential and important concept, especially in system analysis and synthesis. During the past few years, there have been extensive studies. Brian D. O. Anderson, in 1971, had considered the stability of the time-invariant deterministic system in section 3.2 \cite{Anderson1971}. Under some additional assumptions, the stabilization of the closed-loop system was formed when the control law resulting from an infinite-time performance index \cite{Anderson1971}. Subsequently, the stabilization problem of control systems has been extended to stochastic systems, and then extended to time-delay systems, successively; see \cite{H. Zhang2015}-\cite{W. Zhang2004} and references therein. In \cite{H. Zhang2015}, the stabilization of discrete-time systems with delay and multiplicative noise was studied. And it showed that the system was stabilizable in the mean-square sense if and only if an algebraic Riccati-ZXL equation had a particular solution.

The above mentioned results are all belong to classical control theory, that is, the controller has a certain objective to achieve, and the plant to be controlled does not have its own objective. However, in recent years, in the process of socio-economic modeling, it is often necessary to consider the strategic interactions between different decision makers, especially in order to better achieve the objectives of economic planning and policies. This leads to the emergence of the concept and method of stabilization policy theory in the dynamic game literature with hierarchical decision-making structure \cite{M. Johnson}. Moreover, with the popularization of the socialization of network information, it has become particularly important to consider the stabilization of cyber-physical systems for multiply decision makers \cite{M. Zhu2011}. In \cite{M. Zhu2011}, it considered a resilient control problem where a linear dynamic system was subject to cyber attacks launched by correlated jammers in cyber-physical systems. It had proposed a novel leader-follower game formulation to model the interdependency between the operator and adversaries and devised a receding-horizon Stackelberg control law to maintain situational awareness, and further analyzed closed-loop system stability and performance. In \cite{Mukaidania2017}, the infinite horizon linear-quadratic Stackelberg games for a discrete-time stochastic system with multiple decision makers was studied and necessary conditions for the existence of the Stackelberg strategy set were derived in terms of the solvability of cross-coupled stochastic algebraic equations. In
\cite{M. Jungers2008}, the linear-quadratic Stackelberg differential games including time preference rates with an open-loop information structure was investigated. And sufficient conditions to guarantee a predefined degree of stability were given based on the distribution of the eigenvalues in the complex plane. By using a no-memory state feedback representation, \cite{M. Li2002} developed an incentive strategy for discrete time LQ state feedback Stackelberg games. It designed an incentive policy for the leader, while the followers rational reaction guaranteed system stability. \cite{K. Xu2020} proposed an adaptive dynamic programming algorithm by solving the coupled partial differential equations and the Stackelberg feedback equilibrium solution was obtained to ensure the stability of the system according to the Lyapunov function. However, most works are focused on the stabilization of the Stackelberg strategy in infinite horizon. Most recently, \cite{Renren2019} introduced  a new control framework called game-based control systems (GBCSs), which had a hierarchical decision-making structure. And \cite{Renren2019} focused on the linear stochastic systems to give some explicit necessary and sufficient algebraic conditions on the controllability of the Nash equilibrium in GBCSs.

Motivated by \cite{Renren2019}, we study the stabilization problem generated by the leader for the two-player Stackelberg GBCSs, while the follower is designed to minimize its own cost function. In order to address the problem mentioned above, we consider the finite-time horizon open-loop Stackelberg strategy. We optimize the follower firstly, based on the maximum principle, we construct the non-homogeneous relationship between the costate and state and there derives a backward equation. By solving the FBDEs, we derive a homogeneous relationship between the state and non-homogeneous terms in costate equation. In the optimization of the leader, by using the matrix maximum principle, the optimization of the leader is converted into find the optimal gain matrix minimize the cost function and is subject to the state equation.

The rest of the paper is organized as follows. Section II presents some preliminaries of the stabilization problem for the two-player Stackelberg GBCSs. The finite-time horizon results are given in Section III. The stabilization of Stackelberg GBCSs is settled in Section IV. Numerical examples are given in Section V. Conclusions are provided in Section VI.

\emph{Notation:} $R^n$ denotes the $n-$dimensional Euclidean space and $A^T$ denotes the transpose of matrix $A$. A symmetric matrix $M >0$ (reps. $\geq0$) means that it is positive definite (reps. positive semi-definite).

\section{Problem Formulation}
Considering the following discrete-time GBCSs
\begin{eqnarray}\label{opl-eq1}
x_{k+1}=Ax_k+B_1u^1_k+B_2u^2_k,
\end{eqnarray}
where $x_k\in R^n$ is the state, $u^1_k\in R^{m_1}$ and $u^2_k\in R^{m_2}$ are the control inputs of the follower and the leader, respectively. $A$ and $B_i$, $i=1, 2$ are constant matrices of appropriate dimensions. The initial state is $x_0\in R^n$, and the cost function to be minimized by $u^1_k$ is
\begin{eqnarray}
J_1 &\hspace{-0.8em}=&\hspace{-0.8em} \sum\limits_{k=0}^\infty(x^T_kQ_1x_k+u^{1T}_kR_{11}u^1_k+u^{2T}_kR_{12}u^2_k),\label{opl-eq2-1}
\end{eqnarray}
where $Q_1$, $R_{12}$ are positive semi-definite matrices and $R_{11}$ is positive definite matrix of compatible dimensions.

Since (\ref{opl-eq1})-(\ref{opl-eq2-1}) is the linear-quadratic system, we limit the leader's admissible control to the linear feedback strategies, i.e., the admissible control set of the leader is assumed to be $\mathcal{U}_2$, such that
\begin{eqnarray}\label{opl-eq3}
 \mathcal{U}_2=\{u^2_k=K^2x_k: K^2\in R^{m_2\times n}\},
\end{eqnarray}
and the admissible control set of the follower is assumed to be $\mathcal{U}_1$.

Now, we are in the position to give the main problem we want to address.\\
\textbf{Problem 1:} Under the optimal controller  $u^1_k\in \mathcal{U}_1$ which minimizes the cost function $J_1$, find the necessary and
sufficient conditions for $u^2_k=K^2x_k\in \mathcal{U}_2$ to stabilize the GBCS (\ref{opl-eq1}).

What should be noted is that in order to fundamentally solve the Problem 1 mentioned above, we covert the problem of the stabilization of the GBCSs into a optimization problem based on the Stackelberg game, and denote the cost function of the leader as
\begin{eqnarray}
J_2 &\hspace{-0.8em}=&\hspace{-0.8em} \sum\limits_{k=0}^\infty(x^T_kQ_2x_k+u^{1T}_kR_{21}u^1_k+u^{2T}_kR_{22}u^2_k),\label{opl-eq2-2}
\end{eqnarray}
where $Q_2$ is positive semi-definite matrices, such that there exists some matrices $C_2$ satisfying $Q_2=C^T_2C_2$ and $R_{21}$, $R_{22}$ are positive definite matrices of compatible dimensions.

\section{Finite-Horizon Results}

It should pointed out that the controller to minimize (\ref{opl-eq2-2}) will be the best stabilization controller for (\ref{opl-eq1}). Now we will find the controller $u^1_k$ and $u^2_k$ to minimize (\ref{opl-eq1}) and (\ref{opl-eq2-2}), respectively. To this end, we firstly consider the following problem in finite-time horizon.\\
\textbf{Problem 2:} Find the unique open-loop Stackelberg strategy $(u^1_k, u^2_k)$, subject to (\ref{opl-eq1}), where $u^2_k$ is with state feedback form of $u^2_k=K^2_kx_k$, and such that
\begin{eqnarray}\label{opl-eq5}
  J^1_N(u^{1*}(u^2), u^2) &\hspace{-0.8em}\leq&\hspace{-0.8em} J^1_N(u^1(u^2), u^2), \\
  J^2_N(u^{1*}(u^{2*}), u^{2*}) &\hspace{-0.8em}\leq&\hspace{-0.8em}J^2_N(u^{1*}(u^{2}), u^{2}),
\end{eqnarray}
where
\begin{eqnarray}
 J^1_N &\hspace{-0.8em}=&\hspace{-0.8em} \sum\limits_{k=0}^N(x^T_kQ_1x_k+u^{1T}_kR_{11}u^1_k+u^{2T}_kR_{12}u^2_k)+x^T_{N+1}H_1x_{N+1},\label{opl-eq4-1}\\
J^2_N &\hspace{-0.8em}=&\hspace{-0.8em} \sum\limits_{k=0}^N(x^T_kQ_2x_k+u^{1T}_kR_{21}u^1_k+u^{2T}_kR_{22}u^2_k)+x^T_{N+1}H_2x_{N+1}.\label{opl-eq4-2}
\end{eqnarray}

To solve the Problem 2, we introduce the following Riccati equations
\begin{eqnarray}
P^1_k&\hspace{-0.8em}=&\hspace{-0.8em} Q_1+A^TP^1_{k+1}A-A^TP^1_{k+1}B_1(\Gamma^1_{k+1})^{-1}B^{T}_1P^1_{k+1}A,\nonumber\\\label{opl-eq6-1}\\
P^2_k&\hspace{-0.8em}=&\hspace{-0.8em}Q_2+A^TY^T_{k+1}R_{21}Y_{k+1}A+A^TM^{1T}_{k+1}(\Upsilon^{-1}_{k+1})^TP^{2}_{k+1}\Upsilon^{-1}_{k+1}M^1_{k+1}A-M^{2T}_{k+1}(\Gamma^2_{k+1})^{-1} M^2_{k+1},\label{opl-eq6-2}
\end{eqnarray}
with terminal value $P^1_{N+1}=H_1$ and $P^2_{N+1}=H_2$,
where
\begin{eqnarray}
\Gamma^1_{k+1}&\hspace{-0.8em}=&\hspace{-0.8em}R_{11}+B^T_1P^1_{k+1}B_1,\label{opl-eq71}\\
\Gamma^2_{k+1}&\hspace{-0.8em}=&\hspace{-0.8em}R_{22}+B^T_2Y^T_{k+1}R_{21}Y_{k+1}B_2+B^T_2M^{1T}_{k+1}(\Upsilon^{-1}_{k+1})^TP^{2}_{k+1}\Upsilon^{-1}_{k+1}M^1_{k+1}B_2,\label{opl-eq7-2}\\
 M^{1}_{k+1}&\hspace{-0.8em}=&\hspace{-0.8em} I-B_1(\Gamma^1_{k+1})^{-1}B^T_1P^1_{k+1},\label{opl-eq7-3}\\
 M^{2}_{k+1}&\hspace{-0.8em}=&\hspace{-0.8em}B^T_2Y^T_{k+1}R_{21}Y_{k+1}A+B^T_2M^{1T}_{k+1}(\Upsilon^{-1}_{k+1})^TP^{2}_{k+1} \Upsilon^{-1}_{k+1}M^1_{k+1}A,\label{opl-eq7-4}\\
 S_{k+1}&\hspace{-0.8em}=&\hspace{-0.8em}(\Gamma^1_{k+1})^{-1}B^T_1P^1_{k+1},\label{opl-eq7-5}\\
 \Upsilon_{k+1}&\hspace{-0.8em}=&\hspace{-0.8em}I+B_1(\Gamma^1_{k+1})^{-1}B^T_1T_{k+1},\label{opl-eq7-6}\\
 Y_{k+1} &\hspace{-0.8em}=&\hspace{-0.8em}S_{k+1}+(\Gamma^1_{k+1})^{-1}B^T_1T_{k+1}\Upsilon^{-1}_{k+1}M^1_{k+1},\label{opl-eq7-7}\\
 K^1_k&\hspace{-0.8em}=&\hspace{-0.8em}-Y_{k+1}(A+B_2K^2_{k}),\label{opl-eq7-8}\\
 K^2_k&\hspace{-0.8em}=&\hspace{-0.8em}-(\Gamma^2_{k+1})^{-1} M^2_{k+1},\label{opl-eq7-9}\\
 T_k &\hspace{-0.8em}=&\hspace{-0.8em} A^TM^{1T}_{k+1}T_{k+1}\Upsilon^{-1}_{k+1}M^1_{k+1}(A+B_2K^2_{k})+A^TM^{1T}_{k+1}P^1_{k+1}B_2K^2_{k},\label{opl-eq7-10}
\end{eqnarray}
with terminal value $T_{N+1}=0$.

We now give the main results in Theorem 1 to find the open-loop Stackelberg strategy stated in Problem 2.\\
\textbf{Theorem 1:} If $\Upsilon_{k+1}$ is invertible, then the Problem 2 has the unique open-loop Stackelberg strategy. And the optimal Stackelberg strategy for the follower and the leader are
\begin{eqnarray}
u^1_k&\hspace{-0.8em}=&\hspace{-0.8em} K^1_kx_k,\label{opl-eq8-1} \\
 u^2_k &\hspace{-0.8em}=&\hspace{-0.8em} K^2_kx_k.\label{opl-eq8-2}
\end{eqnarray}

The optimal cost functions of the follower and the leader are
\begin{eqnarray}
J^{1*}_N&\hspace{-0.8em}=&\hspace{-0.8em} x^T_0[P^1_0+T^T_0+T_0+\Xi]x_0,\label{opl-eq9-1}\\
J^{2*}_N&\hspace{-0.8em} =&\hspace{-0.8em}x^T_0P^2_0x_0,\label{opl-eq9-2}
\end{eqnarray}
where
\begin{eqnarray}\label{opl-eq10}
\Xi &\hspace{-0.8em}=&\hspace{-0.8em}\sum\limits_{k=0}^N[\Phi^T_{k+1, 1}T^T_{k+1}M^{1}_{k+1}B_2K^2_k \Phi_{k, 1}-\Phi^T_{k+1,1}T^T_{k+1} B_1(\Gamma^1_{k+1})^{-1}B^T_1T_{k+1}\Phi_{k+1, 1}+\Phi^T_{k, 1} K^{2T}_k(R_{12}\nonumber\\
 &\hspace{-0.8em}&\hspace{-0.8em}+B^T_2P^1_{k+1}M^{1}_{k+1}B_2)K^2_k \Phi_{k, 1}
+\Phi^T_{k, 1} K^{2T}_kB^T_2 M^{1T}_{k+1}T_{k+1}\Phi_{k+1, 1}],
\end{eqnarray}
with
\begin{eqnarray*}
 \Phi_{k+1, 1}&\hspace{-0.8em}=&\hspace{-0.8em}A_{k+1}A_k\ldots A_1,\\
 A_{k+1}&\hspace{-0.8em}=&\hspace{-0.8em}\Upsilon^{-1}_{k+1}M^1_{k+1}(A+B_2K^2_{k}).
\end{eqnarray*}

\textbf{Proof:} The proof will be divided into three steps. The first step is to optimize the follower by the maximum principle. And the second step is to optimize the leader by the matrix maximum principle. Finally, we will calculate the optimal cost functions of the two players.

\subsection{The optimization of the follower}
We consider the optimization of the follower. By using the maximum principle, we formulate the non-homogeneous relationship between the state and the costate in this subsection.

To this end, define the Hamiltonian function of the follower as
\begin{eqnarray}\label{opl-eq11}
  H^1_k&\hspace{-0.8em}=&\hspace{-0.8em} x^T_kQ_1x_k+u^{1T}_kR_{11}u^1_k+u^{2T}_kR_{12}u^2_k+\lambda^T_k(Ax_k+B_1u^1_k+B_2u^2_k),
\end{eqnarray}
then according to the maximum principle, we have
\begin{eqnarray}
0&\hspace{-0.8em}=&\hspace{-0.8em}\frac{\partial {H^1_k}}{\partial {u^1_k}}=R_{11}u^1_k+B^T_1\lambda_k,\label{opl-eq12-1}\\
\lambda_{k-1}&\hspace{-0.8em}=&\hspace{-0.8em}\frac{\partial {H^1_k}}{\partial {x_k}}=A^T\lambda_k+Q_1x_k, \label{opl-eq12-2}
\end{eqnarray}
with terminal value $\lambda_{N}=P^1_{N+1}x_{N+1}$.

The existence of $u^2_k$ in (\ref{opl-eq1}) leads to the relationship between $\lambda_{k-1}$ and $x_k$ are no longer homogeneous. Thus, we assume that
\begin{eqnarray}\label{opl-eq13}
 \lambda_{k} &\hspace{-0.8em}=&\hspace{-0.8em}P^1_{k+1}x_{k+1}+\zeta_{k},
\end{eqnarray}
with terminal value $\zeta_{N}=0$.

Adding (\ref{opl-eq13}) into the equilibrium equation (\ref{opl-eq12-1}), then it can be rewritten as
\begin{eqnarray}\label{opl-eq14}
0&\hspace{-0.8em}=&\hspace{-0.8em}R_{11}u^1_k+B^T_1(P^1_{k+1}x_{k+1}+\zeta_{k})\nonumber\\
&\hspace{-0.8em}=&\hspace{-0.8em}(R_{11}+B^T_1P^1_{k+1}B_1)u^1_k+B^T_1P^1_{k+1}Ax_k+B^T_1P^1_{k+1}B_2u^2_k+B^T_1\zeta_{k},
\end{eqnarray}
thus, the control of the follower can be calculated as
\begin{eqnarray}\label{opl-eq15}
  u^1_k&\hspace{-0.8em}=&\hspace{-0.8em}-(\Gamma^1_{k+1})^{-1}(B^T_1P^1_{k+1}Ax_k+B^T_1P^1_{k+1}B_2u^2_k+B^T_1\zeta_{k}).
\end{eqnarray}

Next, we shall proof that (\ref{opl-eq13}) is established for any $0\leq k\leq N$ by the method of inductive hypothesis, where $P^1_k$ satisfies (\ref{opl-eq6-1}) and
$\zeta_{k-1}$ satisfies the following equation,
\begin{eqnarray}
 \zeta_{k-1}&\hspace{-0.8em}=&\hspace{-0.8em}A^TM^{1T}_{k+1}\zeta_{k}+A^TM^{1T}_{k+1}P^1_{k+1}B_2u^2_k,\label{opl-eq16-2}
\end{eqnarray}

According to $\lambda_{N}=P^1_{N+1}x_{N+1}$ and $\zeta_{N}=0$, we can derive that (\ref{opl-eq13}) is established for $k=N$. Given any $s\geq0$, assume that (\ref{opl-eq13}) is established for any $k\geq s$, where $P^1_{k+1}$ and $\zeta_{k}$ satisfy (\ref{opl-eq6-1}) and (\ref{opl-eq16-2}), respectively, we shall show that (\ref{opl-eq13}) also holds for $k=s-1$.

By using (\ref{opl-eq13}) and (\ref{opl-eq12-2}), then it derives
\begin{eqnarray}\label{opl-eq17}
\lambda_{s-1}&\hspace{-0.8em}=&\hspace{-0.8em}A^T(P^1_{s+1}x_{s+1}+\zeta_{s})+Q_1x_s\nonumber\\
&\hspace{-0.8em}=&\hspace{-0.8em}A^TP^1_{s+1}Ax_k+A^TP^1_{s+1}B_1u^1_s+A^TP^1_{s+1}B_2u^2_s+A^T\zeta_{s}+Q_1x_s\nonumber\\
&\hspace{-0.8em}=&\hspace{-0.8em}(Q_1+A^TP^1_{s+1}A-A^TP^1_{s+1}B_1(\Gamma^1_{s+1})^{-1}B^{T}_1P^1_{s+1} A)x_s+A^TM^{1T}_{s+1}\zeta_{s}+A^TM^{1T}_{s+1}P^1_{s+1}B_2u^2_s\nonumber\\
&\hspace{-0.8em}=&\hspace{-0.8em}P^1_sx_s+\zeta_{s-1},
\end{eqnarray}
which indicates that (\ref{opl-eq13}) holds for $k=s-1$ with $P^1_s$ and $\zeta_{s-1}$ satisfying (\ref{opl-eq6-1}) and (\ref{opl-eq16-2}), respectively.

Following from (\ref{opl-eq15}), then the state (\ref{opl-eq1}) can be rewritten as
\begin{eqnarray}\label{opl-eq18}
x_{k+1}=M^{1}_{k+1}Ax_k+M^{1}_{k+1}B_2u^2_k-B_1(\Gamma^1_{k+1})^{-1}B^T_1\zeta_{k}.
\end{eqnarray}

Based on the discussion above, we have yielded the FBDEs (\ref{opl-eq16-2}) and (\ref{opl-eq18}) with initial value $x_0$ and terminal value $\zeta_{N}=0$. Then, the optimization of the leader $u^2_k$ is transformed into the optimization of the FBDEs.

\subsection{The optimization of the leader}

In this subsection, we consider the optimization of the leader. We will iteratively solve the FBDEs firstly. And then, by using the matrix maximum principle, the optimal feedback gain matrix of the leader is yielded.

Adding $u^2_k=K^2_kx_k$ into (\ref{opl-eq16-2}) and (\ref{opl-eq18}), respectively, then there derives the FBDEs
\begin{eqnarray}
x_{k+1}&\hspace{-0.8em}=&\hspace{-0.8em}M^{1}_{k+1}(A+B_2K^2_k)x_k-B_1(\Gamma^1_{k+1})^{-1}B^T_1\zeta_{k}, \label{opl-eq19-1}\\
 \zeta_{k-1}&\hspace{-0.8em}=&\hspace{-0.8em}A^TM^{1T}_{k+1}\zeta_{k}+A^TM^{1T}_{k+1}P^1_{k+1}B_2K^2_kx_k,\label{opl-eq19-2}
\end{eqnarray}
with the initial value $x_0$ and the terminal value $\zeta_{N}=0$.

Then, by using inductive hypothesis, we will proof the FBDEs satisfy the homogeneous relationship
\begin{eqnarray}\label{opl-eq20}
  \zeta_{k-1}=T_{k}x_{k},
\end{eqnarray}
and in this way, we have
\begin{eqnarray}\label{opl-eq21}
x_{k+1}&\hspace{-0.8em}=&\hspace{-0.8em}\Upsilon^{-1}_{k+1}M^1_{k+1}(A+B_2K^2_{k})x_{k}.
\end{eqnarray}

For $k=N$, it yields
\begin{eqnarray*}
x_{N+1}&\hspace{-0.8em}=&\hspace{-0.8em}M^{1}_{N+1}(A+B_2K^2_N)x_N, \nonumber\\
 \zeta_{N-1}&\hspace{-0.8em}=&\hspace{-0.8em}A^TM^{1T}_{N+1}P^1_{N+1}B_2K^2_Nx_N= T_Nx_N,
\end{eqnarray*}
which satisfy (\ref{opl-eq20}) and (\ref{opl-eq21}) for $k=N$.

We take any $n\geq 0$, and assume that $\zeta_{k-1}$ and $x_{k+1}$ are as (\ref{opl-eq20}) and (\ref{opl-eq21}) for all $k\geq n+1$. We show that these conditions will also holds for $k=n$.

For $k=n$, it follows that
\begin{eqnarray*}
x_{n+1}&\hspace{-0.8em}=&\hspace{-0.8em}M^{1}_{n+1}(A+B_2K^2_{n})x_{n}-B_1(\Gamma^1_{n+1})^{-1}B^T_1\zeta_{n}\nonumber\\
&\hspace{-0.8em}=&\hspace{-0.8em}M^{1}_{n+1}(A+B_2K^2_{n})x_{n}-B_1(\Gamma^1_{n+1})^{-1}B^T_1T_{n+1}x_{n+1},
\end{eqnarray*}
then the state $x_{n+1}$ is such that
\begin{eqnarray}\label{opl-eq22}
x_{n+1}=\Upsilon^{-1}_{n+1}M^1_{n+1}(A+B_2K^2_{n})x_{n}£¬
\end{eqnarray}
which is exactly (\ref{opl-eq21}) for $k=n$. And the costate equation $\zeta_{n-1}$ satisfies
\begin{eqnarray}\label{opl-eq23}
 \zeta_{n-1}&\hspace{-0.8em}=&\hspace{-0.8em}A^TM^{1T}_{n+1}\zeta_{n}+A^TM^{1T}_{n+1}P^1_{n+1}B_2K^2_{n}x_{n}\nonumber\\
 &\hspace{-0.8em}=&\hspace{-0.8em}[A^TM^{1T}_{n+1}T_{n+1}\Upsilon^{-1}_{n+1}M^1_{n+1}(A+B_2K^2_{n})+A^TM^{1T}_{n+1}P^1_{n+1}B_2K^2_{n}]x_{n}\nonumber\\
 &\hspace{-0.8em}=&\hspace{-0.8em}T_{n}x_{n},
\end{eqnarray}
which is exactly (\ref{opl-eq20}) for $k=n$.

To this end, the control input of the follower can be written as
\begin{eqnarray}\label{opl-eq24}
 u^1_k&\hspace{-0.8em}=&\hspace{-0.8em}-(\Gamma^1_{k+1})^{-1}(B^T_1P^1_{k+1}Ax_k+B^T_1P^1_{k+1}B_2u^2_k+B^T_1\zeta_{k})\nonumber\\
 &\hspace{-0.8em}=&\hspace{-0.8em}-S_{k+1}(A+B_2K^2_{k})x_k-(\Gamma^1_{k+1})^{-1}B^T_1T_{k+1}\Upsilon^{-1}_{k+1} M^1_{k+1}(A+B_2K^2_{k})x_k\nonumber\\
 &\hspace{-0.8em}=&\hspace{-0.8em} -Y_{k+1}(A+B_2K^2_{k})x_k.
\end{eqnarray}

Now, we are in the position to calculate the optimal feedback gain matrix of the leader.

Adding (\ref{opl-eq24}) and $u^2_k=K^2_kx_k$ into (\ref{opl-eq4-2}), we have
\begin{eqnarray}\label{opl-eq25}
J^2_N &\hspace{-0.8em}=&\hspace{-0.8em}\sum\limits_{k=0}^N[x^T_kQ_2x_k+x^T_k(A+B_2K^2_{k})^TY^T_{k+1}R_{21}Y_{k+1}(A+B_2K^2_{k})x_k+x^T_kK^{2T}_{k}R_{22}K^2_{k}x_k]\nonumber\\
 &\hspace{-0.8em}&\hspace{-0.8em}+x^T_{N+1}P^2_{N+1}x_{N+1}
\end{eqnarray}

Denote the matrix as
\begin{eqnarray}\label{opl-eq26}
  X_{k+1} = x_{k+1} x^T_{k+1}.
\end{eqnarray}
Then, combining (\ref{opl-eq26}) with (\ref{opl-eq21}), it derives that
\begin{eqnarray}\label{opl-eq27}
  X_{k+1} &\hspace{-0.8em}=&\hspace{-0.8em} \Upsilon^{-1}_{k+1}M^1_{k+1}(A+B_2K^2_{k})X_{k} (A+B_2K^2_{k})^T M^{1T}_{k+1}(\Upsilon^{-1}_{k+1})^T.
\end{eqnarray}
Hence, the optimization of the leader can be converted into finding the optimal feedback gain matrix $K^2_k$, minimizes the cost function (\ref{opl-eq25}) and satisfies (\ref{opl-eq27}).

To this end, denote the Hamiltonian function of the leader as
\begin{eqnarray}\label{opl-eq28}
  H^2_k &\hspace{-0.8em}=&\hspace{-0.8em}Tr\Big([ Q_2+(A+B_2K^2_{k})^TY^T_{k+1}R_{21}Y_{k+1}(A+B_2K^2_{k})+K^{2T}_{k}R_{22}K^2_{k}]X_k+\Upsilon^{-1}_{k+1}M^1_{k+1}(A+B_2K^2_{k}) \nonumber\\
 &\hspace{-0.8em}&\hspace{-0.8em}\times X_{k}(A+B_2K^2_{k})^TM^{1T}_{k+1}(\Upsilon^{-1}_{k+1})^TP^{2T}_{k+1}\Big),
\end{eqnarray}
applying the matrix maximum principle, it yields
\begin{eqnarray}
  0&\hspace{-0.8em}=&\hspace{-0.8em}\frac{\partial {H^2_k}}{\partial {K^2_k}} \nonumber\\
 &\hspace{-0.8em}=&\hspace{-0.8em}R_{22}K^2_kX_k+R_{22}K^2_kX^T_k+B^T_2Y^T_{k+1}R_{21}Y_{k+1}AX^T_k+B^T_2Y^T_{k+1}R_{21}Y_{k+1}AX_k\nonumber\\
 &\hspace{-0.8em}&\hspace{-0.8em}+B^T_2Y^T_{k+1}R_{21}Y_{k+1}B_2K^2_kX_k+B^T_2Y^T_{k+1}R_{21}Y_{k+1}B_2K^2_kX^T_k+B^T_2M^{1T}_{k+1}(\Upsilon^{-1}_{k+1})^T\nonumber\\
 &\hspace{-0.8em}&\hspace{-0.8em}\times P^{2T}_{k+1} \Upsilon^{-1}_{k+1}M^1_{k+1}AX_k
 +B^T_2M^{1T}_{k+1}(\Upsilon^{-1}_{k+1})^T P^{2}_{k+1}\Upsilon^{-1}_{k+1}M^1_{k+1}AX^T_k
 +B^T_2M^{1T}_{k+1}(\Upsilon^{-1}_{k+1})^T\nonumber\\
 &\hspace{-0.8em}&\hspace{-0.8em}\times P^{2T}_{k+1}\Upsilon^{-1}_{k+1}M^1_{k+1}B_2K^2_kX_k
 +B^T_2M^{1T}_{k+1}(\Upsilon^{-1}_{k+1})^TP^{2}_{k+1} \Upsilon^{-1}_{k+1}M^1_{k+1}B_2K^2_kX^T_k,\label{opl-eq29-1}\\
 P^2_k&\hspace{-0.8em}=&\hspace{-0.8em}\frac{\partial {H^2_k}}{\partial {X_k}} \nonumber\\
 &\hspace{-0.8em}=&\hspace{-0.8em}Q_2+(A+B_2K^2_{k})^TY^T_{k+1}R_{21}Y_{k+1}(A+B_2K^2_{k})+K^{2T}_{k}R_{22}K^2_{k}+(A+B_2K^2_{k})^TM^{1T}_{k+1}(\Upsilon^{-1}_{k+1})^T\nonumber\\
 &\hspace{-0.8em}&\hspace{-0.8em}\times P^{2}_{k+1}\Upsilon^{-1}_{k+1}M^1_{k+1}(A+B_2K^2_{k}),\label{opl-eq29-2}
\end{eqnarray}
with terminal value $P^2_{N+1}=H_2$.

According to the symmetric of the matrices $X_k$ and $P^{2}_{k+1}$, we have,
\begin{eqnarray}\label{opl-eq30}
 0&\hspace{-0.8em}=&\hspace{-0.8em}[R_{22}+B^T_2Y^T_{k+1}R_{21}Y_{k+1}B_2+B^T_2M^{1T}_{k+1}(\Upsilon^{-1}_{k+1})^TP^{2}_{k+1}\Upsilon^{-1}_{k+1}M^1_{k+1}B_2]K^2_k\nonumber\\
 &\hspace{-0.8em}&\hspace{-0.8em} +B^T_2Y^T_{k+1}R_{21}Y_{k+1}A+B^T_2M^{1T}_{k+1}(\Upsilon^{-1}_{k+1})^TP^{2}_{k+1}\Upsilon^{-1}_{k+1}M^1_{k+1}A\nonumber\\
&\hspace{-0.8em}=&\hspace{-0.8em}\Gamma^2_{k+1}K^2_k+M^2_{k+1}.
\end{eqnarray}

Since $\Gamma^2_{k+1}>0$, then the optimal gain matrix of the leader can be uniquely obtained, which satisfies (\ref{opl-eq7-9}).
Based on the discussion above, we can derive that the optimal controllers of the follower and the leader are exactly (\ref{opl-eq8-1}) and (\ref{opl-eq8-2}).

By using (\ref{opl-eq7-9}), the Riccait equation of the leader can be written as
\begin{eqnarray}\label{opl-eq33}
P^2_k&\hspace{-0.8em}=&\hspace{-0.8em}Q_2+K^{2T}_{k}\Gamma^2_{k+1}K^2_{k}+K^{2T}_{k}M^2_{k+1}+M^{2T}_{k+1}K^2_{k}+A^TY^T_{k+1}R_{21}Y_{k+1}A\nonumber\\
 &\hspace{-0.8em}&\hspace{-0.8em}+A^TM^{1T}_{k+1}(\Upsilon^{-1}_{k+1})^TP^{2}_{k+1}\Upsilon^{-1}_{k+1} M^1_{k+1}A\nonumber\\
&\hspace{-0.8em}=&\hspace{-0.8em}Q_2+A^TY^T_{k+1}R_{21}Y_{k+1}A+A^TM^{1T}_{k+1}(\Upsilon^{-1}_{k+1})^TP^{2}_{k+1}\Upsilon^{-1}_{k+1}M^1_{k+1}A-M^{2T}_{k+1}(\Gamma^2_{k+1})^{-1} M^2_{k+1},
\end{eqnarray}
which is exactly (\ref{opl-eq6-2}), with terminal value $P^{2}_{N+1}=H_2$.

Lastly, we will calculate the optimal cost functions of the two players. We shall calculate the optimal cost function of the follower firstly.

By applying (\ref{opl-eq1}) and (\ref{opl-eq12-2}), we get
\begin{eqnarray*}
  x^T_k\lambda_{k-1}-x^T_{k+1}\lambda_{k} &\hspace{-0.8em}=&\hspace{-0.8em}x^T_k(A^T\lambda_k+Q_1x_k)-(Ax_k+B_1u^1_k+B_2u^2_k)^T\lambda_{k}\nonumber\\
&\hspace{-0.8em}=&\hspace{-0.8em}x^T_kQ_1x_k-(B_1u^1_k+B_2u^2_k)^T\lambda_{k},
\end{eqnarray*}
Adding the above equation from $k=0$ to $k=N$, we have
\begin{eqnarray}\label{opl-eq34}
x^T_0\lambda_{-1}-x^T_{N+1}\lambda_{N} &\hspace{-0.8em}=&\hspace{-0.8em}\sum\limits_{k=0}^N[x^T_kQ_1x_k-(B_1u^1_k+B_2u^2_k)^T\lambda_{k}].
\end{eqnarray}

Compared (\ref{opl-eq34}) with the cost function (\ref{opl-eq2-1}) and based on the equilibrium condition (\ref{opl-eq12-1}), there derives
\begin{eqnarray}\label{opl-eq35}
  J^1_N &\hspace{-0.8em}=&\hspace{-0.8em} x^T_0\lambda_{-1}+ \sum\limits_{k=0}^N[u^{1T}_k(R_{11}u^1_k+B^T_1\lambda_k)+u^{2T}_k(R_{12}u^2_k+B^T_2\lambda_k)]\nonumber\\
&\hspace{-0.8em}=&\hspace{-0.8em}x^T_0\lambda_{-1}+ \sum\limits_{k=0}^N[u^{2T}_k(R_{12}u^2_k+B^T_2\lambda_k)].
\end{eqnarray}
And by using (\ref{opl-eq13}), (\ref{opl-eq16-2}) and (\ref{opl-eq18}), it yields
\begin{eqnarray}\label{opl-eq36}
 u^{2T}_kB^T_2\lambda_k &\hspace{-0.8em}=&\hspace{-0.8em} u^{2T}_kB^T_2(P^1_{k+1}x_{k+1}+\zeta_k)\nonumber\\
&\hspace{-0.8em}=&\hspace{-0.8em}u^{2T}_kB^T_2P^1_{k+1}(M^{1}_{k+1}Ax_k+M^{1}_{k+1}B_2u^2_k-B_1(\Gamma^1_{k+1})^{-1}B^T_1\zeta_{k})+u^{2T}_kB^T_2\zeta_k\nonumber\\
&\hspace{-0.8em}=&\hspace{-0.8em}(\zeta_{k-1}-A^TM^{1T}_{k+1}\zeta_{k})^Tx_k+u^{2T}_kB^T_2P^1_{k+1} M^{1}_{k+1}B_2u^2_k+u^{2T}_kB^T_2M^{1T}_{k+1}\zeta_k\nonumber\\ &\hspace{-0.8em}=&\hspace{-0.8em}\zeta^T_{k-1}x_k-\zeta^T_{k}[x_{k+1}-M^{1}_{k+1}B_2u^2_k+B_1 (\Gamma^1_{k+1})^{-1}B^T_1\zeta_{k}]+u^{2T}_kB^T_2P^1_{k+1}M^{1}_{k+1}B_2u^2_k\nonumber\\
&\hspace{-0.8em}&\hspace{-0.8em}+u^{2T}_kB^T_2M^{1T}_{k+1}\zeta_k\nonumber\\
&\hspace{-0.8em}=&\hspace{-0.8em}\zeta^T_{k-1}x_k-\zeta^T_{k}x_{k+1}+\zeta^T_{k}M^{1}_{k+1}B_2u^2_k-\zeta^T_{k}B_1(\Gamma^1_{k+1})^{-1}B^T_1\zeta_{k}+u^{2T}_kB^T_2P^1_{k+1}M^{1}_{k+1}B_2u^2_k\nonumber\\
&\hspace{-0.8em}&\hspace{-0.8em}+u^{2T}_kB^T_2M^{1T}_{k+1}\zeta_k.
\end{eqnarray}

Combining (\ref{opl-eq35}) with (\ref{opl-eq36}), thus, we have
\begin{eqnarray}\label{opl-eq37}
J^1_N &\hspace{-0.8em}=&\hspace{-0.8em} x^T_0\lambda_{-1}+ \zeta^T_{-1}x_0+\sum\limits_{k=0}^N[\zeta^T_{k}M^{1}_{k+1}B_2u^2_k-\zeta^T_{k}B_1(\Gamma^1_{k+1})^{-1}B^T_1\zeta_{k}\nonumber\\
&\hspace{-0.8em}&\hspace{-0.8em}+u^{2T}_k(R_{22}+B^T_2P^1_{k+1}M^{1}_{k+1}B_2)u^2_k+u^{2T}_kB^T_2M^{1T}_{k+1}\zeta_k].
\end{eqnarray}

According to (\ref{opl-eq21}), we can yield that
\begin{eqnarray}\label{opl-eq38}
  x_{k+1}&\hspace{-0.8em}=&\hspace{-0.8em}\Upsilon^{-1}_{k+1}M^1_{k+1}(A+B_2K^2_{k})x_{k}= A_{k+1}x_k\nonumber\\
&\hspace{-0.8em}=&\hspace{-0.8em}A_{k+1}A_kx_{k-1}=\ldots=A_{k+1}A_k\ldots A_1x_0=\Phi_{k+1, 1}x_0.
\end{eqnarray}
Adding it into (\ref{opl-eq20}) and $u^2_k=K^2_kx_k$, there follows
\begin{eqnarray}
 \zeta_k &\hspace{-0.8em}=&\hspace{-0.8em} T_{k+1}\Phi_{k+1, 1}x_0, \\
 u^2_k &\hspace{-0.8em}=&\hspace{-0.8em} K^2_k \Phi_{k, 1}x_0.
\end{eqnarray}

Then the optimal cost function of the follower is
\begin{eqnarray}\label{opl-eq40}
J^{1*}_N &\hspace{-0.8em}=&\hspace{-0.8em} x^T_0P^1_0x_0+x^T_0\zeta_{-1}+ \zeta^T_{-1}x_0+\sum\limits_{k=0}^N\Big(x^T_0[\Phi^T_{k+1, 1}T^T_{k+1} M^{1}_{k+1}B_2K^2_k \Phi_{k, 1}-\Phi^T_{k+1,1}T^T_{k+1}B_1\nonumber\\
&\hspace{-0.8em}&\hspace{-0.8em}\times (\Gamma^1_{k+1})^{-1}B^T_1T_{k+1}\Phi_{k+1, 1}
+\Phi^T_{k, 1} K^{2T}_k(R_{22}+B^T_2P^1_{k+1}M^{1}_{k+1}B_2)K^2_k \Phi_{k, 1}\nonumber\\
&\hspace{-0.8em}&\hspace{-0.8em}+\Phi^T_{k, 1} K^{2T}_k B^T_2M^{1T}_{k+1}T_{k+1}\Phi_{k+1, 1}]x_0\Big)\nonumber\\
&\hspace{-0.8em}=&\hspace{-0.8em}x^T_0[P^1_0+T^T_0+T_0+\Xi]x_0.
\end{eqnarray}

Finally, we will consider the optimal cost function of the leader.

Adding
\begin{eqnarray}\label{opl-eq41}
  u^1_k&\hspace{-0.8em}=&\hspace{-0.8em} -Y_{k+1}(A+B_2K^2_{k})x_k\nonumber\\
&\hspace{-0.8em}=&\hspace{-0.8em}-Y_{k+1}Ax_k-Y_{k+1}B_2u^2_{k}
\end{eqnarray}
into the state (\ref{opl-eq1}), we have
\begin{eqnarray}\label{opl-eq42}
  x_{k+1}&\hspace{-0.8em}=&\hspace{-0.8em} (A-B_1Y_{k+1}A)x_k+(B_2-B_1Y_{k+1}B_2)u^2_{k}.
\end{eqnarray}

Denote the value function of the leader as
\begin{eqnarray}\label{opl-eq43}
V^2_k &\hspace{-0.8em}=&\hspace{-0.8em}x^T_{k+1}P^2_{k+1}x_{k+1}+x^T_kQ_2x_k+u^{1T}_kR_{21}u^1_k+u^{2T}_kR_{22}u^2_k.
\end{eqnarray}

Combining (\ref{opl-eq42}) with (\ref{opl-eq43}), the value function can be written as
\begin{eqnarray}\label{opl-eq44}
V^2_k&\hspace{-0.8em}=&\hspace{-0.8em} u^{2T}_k[R_{22}+B^T_2Y^T_{k+1}R_{21}Y_{k+1}B_2+B^T_2P^2_{k+1}B_2-B^T_2P^2_{k+1}B_1Y_{k+1}B_2-B^T_2Y^T_{k+1}B^T_1\nonumber\\
&\hspace{-0.8em}&\hspace{-0.8em}\times P^2_{k+1}B_2+B^T_2Y^T_{k+1}B^T_1P^2_{k+1}B_1Y_{k+1}B_2]u^2_k+2u^{2T}_k[B^T_2P^2_{k+1}A-B^T_2P^2_{k+1}B_1Y_{k+1}A\nonumber\\
&\hspace{-0.8em}&\hspace{-0.8em}-B^T_2Y^T_{k+1}B^T_1P^2_{k+1}A+B^T_2Y^T_{k+1}B^T_1P^2_{k+1}B_1Y_{k+1}A+B^T_2Y^T_{k+1}R_{21}Y_{k+1}A]x_k\nonumber\\
&\hspace{-0.8em}&\hspace{-0.8em}+x^T_k[Q_2+A^TP^2_{k+1}A-A^TP^2_{k+1}B_1Y_{k+1}A-A^TY^T_{k+1}B^T_1P^2_{k+1}A+A^TY^T_{k+1}B^T_1\nonumber\\
&\hspace{-0.8em}&\hspace{-0.8em}\times P^2_{k+1}B_1Y_{k+1}A+A^TY^T_{k+1}R_{21}Y_{k+1}A]x_k\nonumber\\
&\hspace{-0.8em}=&\hspace{-0.8em}u^{2T}_k[R_{22}+B^T_2Y^T_{k+1}R_{21}Y_{k+1}B_2+B^T_2(I-B_1Y_{k+1})^TP^2_{k+1}(I-B_1Y_{k+1})B_2]u^{2}_k\nonumber\\
&\hspace{-0.8em}&\hspace{-0.8em}+2u^{2T}_k[B^T_2(I-B_1Y_{k+1})^TP^2_{k+1}(I-B_1Y_{k+1})A+B^T_2Y^T_{k+1}R_{21}Y_{k+1}A]x_k\nonumber\\
&\hspace{-0.8em}&\hspace{-0.8em}+x^T_k[Q_2+A^T(I-B_1Y_{k+1})^TP^2_{k+1}(I-B_1Y_{k+1})A+A^TY^T_{k+1}R_{21}Y_{k+1}A]x_k\nonumber\\
&\hspace{-0.8em}=&\hspace{-0.8em}[u^{2}_k+(\Gamma^2_{k+1})^{-1}M^2_{k+1}x_k]^T\Gamma^2_{k+1}[u^{2}_k+(\Gamma^2_{k+1})^{-1}M^2_{k+1}x_k]
+x^T_kP^2_kx_k,
\end{eqnarray}
where the last equality is established because of
\begin{eqnarray*}
I-B_1Y_{k+1}&\hspace{-0.8em}=&\hspace{-0.8em}M^1_{k+1}-(I+B_1(\Gamma^1_{k+1})^{-1}B^T_1T_{k+1})\Upsilon^{-1}_{k+1} M^1_{k+1}+\Upsilon^{-1}_{k+1}M^1_{k+1}\nonumber\\
&\hspace{-0.8em}=&\hspace{-0.8em}\Upsilon^{-1}_{k+1}M^1_{k+1}.
\end{eqnarray*}

Adding (\ref{opl-eq44}) from $k=0$ to $k=N$ on both sides of the equation, it yields that
\begin{eqnarray}\label{opl-eq45}
&\hspace{-0.8em}&\hspace{-0.8em}x^T_{N+1}P^2_{N+1}x_{N+1}-x^T_0P^2_0x_0+\sum\limits_{k=0}^N(x^T_kQ_2x_k+u^{1T}_kR_{21}u^1_k+u^{2T}_kR_{22}u^2_k)\nonumber\\
&\hspace{-0.8em}=&\hspace{-0.8em}\sum\limits_{k=0}^N[u^{2}_k+(\Gamma^2_{k+1})^{-1}M^2_{k+1}x_k]^T\Gamma^2_{k+1}[u^{2}_k+(\Gamma^2_{k+1})^{-1} M^2_{k+1}x_k],
\end{eqnarray}
according to the optimal control of the leader is $u^{2}_k=-(\Gamma^2_{k+1})^{-1}M^2_{k+1}x_k$, and compared with (\ref{opl-eq2-2}), the optimal cost function of the leader is exactly (\ref{opl-eq9-2}).
This completest the proof of Theorem 1.
$\hfill\blacksquare$

\section{Stabilization Results}

In this section, the stabilization of Stackelberg GBCSs will be investigated.

Before giving the main results, we introduce the following AREs:
\begin{eqnarray}
P^1&\hspace{-0.8em}=&\hspace{-0.8em} Q_1+A^TP^1A-A^TP^1B_1(\Gamma^1)^{-1}B^{T}_1P^1A,\label{opl-eq47-1}\\
P^2&\hspace{-0.8em}=&\hspace{-0.8em}Q_2+A^TY^TR_{21}YA+A^TM^{1T}(\Upsilon^{-1})^TP^{2}\Upsilon^{-1}M^1A-M^{2T}(\Gamma^2)^{-1} M^2.\label{opl-eq47-2}
\end{eqnarray}
where
\begin{eqnarray}
\Gamma^1&\hspace{-0.8em}=&\hspace{-0.8em}R_{11}+B^T_1P^1B_1,\label{opl-eq50-1}\\
\Gamma^2&\hspace{-0.8em}=&\hspace{-0.8em}R_{22}+B^T_2Y^TR_{21}YB_2+B^T_2M^{1T}(\Upsilon^{-1})^TP^{2}\Upsilon^{-1}M^1B_2,\label{opl-eq50-2}\\
 M^{1}&\hspace{-0.8em}=&\hspace{-0.8em} I-B_1(\Gamma^1)^{-1}B^T_1P^1,\label{opl-eq50-3}\\
 M^{2}&\hspace{-0.8em}=&\hspace{-0.8em}B^T_2Y^TR_{21}YA+B^T_2M^{1T}(\Upsilon^{-1})^TP^{2}\Upsilon^{-1}M^1A,\label{opl-eq50-4}\\
 S&\hspace{-0.8em}=&\hspace{-0.8em}(\Gamma^1)^{-1}B^T_1P^1,\label{opl-eq50-5}\\
 \Upsilon&\hspace{-0.8em}=&\hspace{-0.8em}I+B_1(\Gamma^1)^{-1}B^T_1T,\label{opl-eq50-6}\\
 Y &\hspace{-0.8em}=&\hspace{-0.8em}S+(\Gamma^1)^{-1}B^T_1T\Upsilon^{-1}M^1,\label{opl-eq50-7}\\
 K^1&\hspace{-0.8em}=&\hspace{-0.8em}-Y(A+B_2K^2),\label{opl-eq50-8}\\
 K^2&\hspace{-0.8em}=&\hspace{-0.8em}-(\Gamma^2)^{-1} M^2,\label{opl-eq50-9}\\
 T &\hspace{-0.8em}=&\hspace{-0.8em} A^TM^{1T}T\Upsilon^{-1}M^1(A+B_2K^2)+A^TM^{1T}P^1B_2K^2.\label{opl-eq50-10}
\end{eqnarray}

Throughout the rest of this paper, the following two standard assumptions are made, readers may refer to \cite{Y. Huang2006}.\\
\textbf{Assumption 1:}  $(A, B_1)$ is controllable.\\
\textbf{Assumption 2:}  $(A, Q^{\frac{1}{2}}_2)$ is observable.

Now, we will give the main results.\\
\textbf{Theorem 2:} Under Assumptions 1, 2 and if $\Upsilon$ is invertible, $u^2_k \in \mathcal{U}_2$ stabilizes the GBCS (\ref{opl-eq1}) in stabilization if and only if there exists a unique positive definite solution $P^1$ and $P^2$ to the AREs (\ref{opl-eq47-1})-(\ref{opl-eq47-2}).
In this case, the stabilizing controller is given as
\begin{eqnarray}
 u^2_k &\hspace{-0.8em}=&\hspace{-0.8em} K^2x_k,\label{opl-eq48}
\end{eqnarray}
while the optimal controller of minimizing $J_1$ is given as following, i.e.,
\begin{eqnarray}
u^1_k&\hspace{-0.8em}=&\hspace{-0.8em} K^1x_k, \label{opl-eq49}
\end{eqnarray}

And the optimal cost functions are
\begin{eqnarray}
 J^*_1&\hspace{-0.8em}=&\hspace{-0.8em}x^T_0(P^1+T+T^T+\Xi)x_0,\label{opl-eq51}\\
 J^*_2&\hspace{-0.8em}=&\hspace{-0.8em}x^T_0P^2x_0,\label{opl-eq51-1},
\end{eqnarray}
where
\begin{eqnarray}
 \Xi&\hspace{-0.8em}=&\hspace{-0.8em}\sum\limits_{k=0}^\infty [(\bar{A}^k)^TK^{2T}(R_{12}+B^T_2P^1M^1B_2)K^2\bar{A}^k+(\bar{A}^k)^TK^{2T}B^T_2M^{1T}T\bar{A}^{k+1}\nonumber\\
 &\hspace{-0.8em}&\hspace{-0.8em}+(\bar{A}^{k+1})^TT^TM^1B_2K^2\bar{A}^k-(\bar{A}^{k+1})^TT^TB_1(\Gamma^1)^{-1}B^T_1T\bar{A}^{k+1}],
\end{eqnarray}
with
\begin{eqnarray}
 \bar{A}&\hspace{-0.8em}=&\hspace{-0.8em}\gamma^{-1}M^1(A+B_2K^2).
\end{eqnarray}

\textbf{Proof:}  We give the necessity proof at first.

\emph{Necessity:} Under Assumption 1 and 2, and if $\Upsilon$ is invertible, suppose there exists controllers $u^1_k$ and $u^2_k$ with constants $K^1$ and $K^2$ such that $u^1_k$ minimizes the cost function $J_1$ and $u^2_k$ stabilizes the system (\ref{opl-eq1}) in stabilization. We will show that there exists a unique solution $P^1>0$ and $P^2>0$ are the solutions to the AREs (\ref{opl-eq47-1}) and (\ref{opl-eq47-2}), respectively.

To make the time horizon $N$ explicit in the finite horizon case, we rewrite $P^1_k$, $P^2_k$, $\Gamma^1_{k+1}$, $\Gamma^2_{k+1}$, $M^{1}_{k+1}$, $M^{2}_{k+1}$, $Y_{k+1}$, $S_{k+1}$, $\Upsilon_{k+1}$, $K^1_k$, $K^2_k$ and $T_k$ in (\ref{opl-eq6-1})-(\ref{opl-eq7-10}) as $P^1_k(N)$, $P^2_k(N)$, $\Gamma^1_{k+1}(N)$, $\Gamma^2_{k+1}(N)$, $M^{1}_{k+1}(N)$, $M^{2}_{k+1}(N)$, $Y_{k+1}(N)$, $S_{k+1}(N)$, $\Upsilon_{k+1}(N)$, $K^1_k(N)$, $K^2_k(N)$ and $T_k(N)$. To consider the infinite horizon case, the terminal weighting matrices in (\ref{opl-eq4-1}) and (\ref{opl-eq4-2}) are set to be zero, i.e., $H_1=H_2=0$.

$\Upsilon_{k+1}(N)$ is invertible and it is noted from Theorem 1, we can conclude that Problem 2 admits a unique solution.

Since $P^1_k(N)$ satisfies the standard Riccati equation
\begin{eqnarray*}
P^1_k(N)&\hspace{-0.8em}=&\hspace{-0.8em} Q_1+A^TP^1_{k+1}(N)A-A^TP^1_{k+1}(N)B_1(\Gamma^1_{k+1}(N))^{-1}B^{T}_1P^1_{k+1}(N)A,
\end{eqnarray*}
under the Assumption 1, and based on the monotone bounded theorem, we can yield that
\begin{eqnarray}
  \lim_{N\rightarrow\infty} P^1_k(N)=P^1>0.
\end{eqnarray}

Next, we will show $P^2>0$.

Rewritten the Riccati equation (\ref{opl-eq6-2}) as
\begin{eqnarray}\label{opl-eq52}
P^2_k(N)&\hspace{-0.8em}=&\hspace{-0.8em}Q_2+A^TY^T_{k+1}(N)R _{21}Y_{k+1}(N)A+A^TM^{1T}_{k+1}(N)(\Upsilon^{-1}_{k+1}(N))^TP^{2}_{k+1}(N)\nonumber\\
&\hspace{-0.8em}&\hspace{-0.8em}\times \Upsilon^{-1}_{k+1}(N)M^1_{k+1}(N)A-M^{2T}_{k+1}(N)(\Gamma^2_{k+1}(N))^{-1} M^2_{k+1}(N)\nonumber\\
&\hspace{-0.8em}=&\hspace{-0.8em}Q_2+A^TY^T_{k+1}(N)R_{21}Y_{k+1}(N)A+A^TM^{1T}_{k+1}(N)(\Upsilon^{-1}_{k+1}(N))^TP^{2}_{k+1}(N)\Upsilon^{-1}_{k+1}(N)\nonumber\\
&\hspace{-0.8em}&\hspace{-0.8em}\times M^1_{k+1}(N)A+K^{2T}_k(N)\Gamma^2_{k+1}(N)K^{2}_k(N)-M^{2T}_{k+1}(N)K^{2}_k(N)-K^{2T}_k(N)M^{2}_{k+1}(N)\nonumber\\
&\hspace{-0.8em}=&\hspace{-0.8em}Q_2+K^{2T}_k(N)R_{22}K^{2}_k(N)+[Y_{k+1}(N)A-Y_{k+1}(N)B_2K^{2}_k(N)]^TR_{21}[Y_{k+1}(N)A\nonumber\\
&\hspace{-0.8em}&\hspace{-0.8em}-Y_{k+1}(N)B_2K^{2}_k(N)]+[\Upsilon^{-1}_{k+1}(N)M^1_{k+1}(N)A-\Upsilon^{-1}_{k+1}(N)M^1_{k+1}(N)B_2K^{2}_k(N)]^T\nonumber\\
&\hspace{-0.8em}&\hspace{-0.8em}\times P^{2}_{k+1}(N)[\Upsilon^{-1}_{k+1}(N)M^1_{k+1}(N)A-\Upsilon^{-1}_{k+1}(N)M^1_{k+1}(N)B_2K^{2}_k(N)].
\end{eqnarray}
According to the terminal condition $H_2= 0$, then by using induction method, we have $P^2_k(N) \geq 0$.

Following from (\ref{opl-eq9-2}), it derives that,
\begin{eqnarray}\label{opl-eq53}
x^T_0P^2_0(N)x_0&\hspace{-0.8em}=&\hspace{-0.8em}J^*_2(N)J^*_2(N+1)= x^T_0P^2_0(N+1)x_0,
\end{eqnarray}
the arbitrary of $x_0$ implies that $P^2_0(N)\leq P^2_0(N+1)$, i.e., $P^2_0(N)$ increases with respect to $N$. Next the boundedness of $P^2_0(N)$ is to be proved as below.

Since the controllers $u^2_k=K^2x_k$ stabilizes system (\ref{opl-eq1}), while $u^1_k=K^1x_k$ minimizes the cost function (\ref{opl-eq2-1}), thus we have
\begin{eqnarray}\label{opl-eq53-1}
 \lim_{k\rightarrow\infty}x_{k} = \lim_{k\rightarrow\infty}(A+B_1K^1+B_2K^2)x_{k-1}=0,
\end{eqnarray}
and then we can yield that
\begin{eqnarray}\label{opl-eq54}
\lim_{k\rightarrow\infty}x^T_kx_k=0.
\end{eqnarray}

Following from \cite{Bouhtouri1999}, it follows that there exist constant $c_1$ satisfying:
\begin{eqnarray*}
  \sum\limits_{k=0}^\infty x^T_kx_k\leq c_1x^T_0x_0.
\end{eqnarray*}
Therefore, noting from (\ref{opl-eq2-2}) that $Q_2$, $K^{1T}R_{21}K^1$ and $K^{2T}R_{22}K^2$ are both bounded, there exists constant $c_2$ such that
\begin{eqnarray}\label{opl-eq55}
J_2 &\hspace{-0.8em}=&\hspace{-0.8em} \sum\limits_{k=0}^\infty(x^T_kQ_2x_k+u^{1T}_kR_{21}u^1_k+u^{2T}_kR_{22}u^2_k)\nonumber\\
&\hspace{-0.8em}=&\hspace{-0.8em}\sum\limits_{k=0}^\infty[x^T_k(Q_2+K^{1T}R_{21}K^1+K^{2T}R_{22}K^2)x_k]\nonumber\\
&\hspace{-0.8em}\leq&\hspace{-0.8em}c_2x^T_0x_0.
\end{eqnarray}

Thus for any $N > 0$, from Theorem 1 we know that
\begin{eqnarray*}
 x^T_0P^2_0(N)x_0=J^*_2(N)\leq J_2\leq c_2x^T_0x_0
\end{eqnarray*}
which indicates that $P^2_0(N)$ is bounded.

Also we know that $P^2_0(N)$ is monotonically increasing, hence $P^2_0(N)$ is convergent, that is
\begin{eqnarray}\label{opl-eq56}
 \lim_{N\rightarrow\infty}P^2_0(N)=P^2.
\end{eqnarray}
According to $P^2_k(N)=P^2_0(N-k)$, we have
\begin{eqnarray}\label{opl-eq57}
\lim_{N\rightarrow\infty}P^2_k(N)=\lim_{N\rightarrow\infty}P^2_0(N-k)=P^2.
\end{eqnarray}

What should be noted is that we do not use the convergence of the $T_k(N)$ in the proof of the convergence of $P^2_k(N)$, however, according to the Riccati equation (\ref{opl-eq6-2}), the existence of the convergence for $P^2_k(N)$ means that $T_k(N)$ is also convergence, thus, we assume $\lim_{N\rightarrow\infty}T_k(N)=T$.

Taking limitation of $N$ on both sides of (\ref{opl-eq6-1})-(\ref{opl-eq7-10}), therefore, (\ref{opl-eq50-1})-(\ref{opl-eq50-10}) can be, respectively, obtained. And $P^1$, $P^2$ satisfy the Riccati equations (\ref{opl-eq47-1}) and (\ref{opl-eq47-2}).

Next, we will show that there exist a positive integer $N_0$ such that $P^2_0(N_0)$ is positive definite. Suppose this is not the case. Since $P^2_0(N)\geq0$, then we can get an non-empty set
\begin{eqnarray}\label{opl-eq58}
  X_N =\{x\in R^n: x\neq0, x^TP^2_0(N_0)x=0\}.
\end{eqnarray}
By using (\ref{opl-eq53}), which means that
\begin{eqnarray*}
 x^T_0P^2_0(N+1)x_0=0, \Rightarrow x^T_0P^2_0(N)x_0=0,
\end{eqnarray*}
i.e., $X_{N+1}\subseteq X_N$. Each $X_N$ is a non-empty finite-dimensional set, so
\begin{eqnarray*}
  1\leq\cdots \leq dim(X_1)\leq dim(X_0)\leq n
\end{eqnarray*}
where dim represents the dimension of the set. Thus, there must exist $N_1$, such that for any $N\geq N_1$,
\begin{eqnarray*}
  dim(X_{N})= dim(X_{N_1})
\end{eqnarray*}
which yields that $X_{N}=X_{N_1}$, and thus
\begin{eqnarray*}
  \bigcap_{N\geq0} X_N= X_{N_1}\neq0.
\end{eqnarray*}
So there exists a nonzero vector $x\in X_{N_1}$, such that
\begin{eqnarray*}
x^TP^2_0(N)x=0, \forall N\geq 0.
\end{eqnarray*}

Let $x_0$ be equal to $x$. Then the optimal value of (\ref{opl-eq2-2}) is as
\begin{eqnarray*}
J^*_2(N) &\hspace{-0.8em}=&\hspace{-0.8em} \sum\limits_{k=0}^N[x^{*T}_kQ_2x^*_k+(u^{1*}_k)^TR_{21}u^{1*}_k+(u^{2*}_k)^TR_{22}u^{2*}_k]\nonumber\\
&\hspace{-0.8em}=&\hspace{-0.8em}x^TP^2_0(N)x=0
\end{eqnarray*}
where $x^*_k$, $u^{1*}_k$ and $u^{2*}_k$ represent the optimal state and the optimal controllers, respectively. Note that $R_{21}>0$, $R_{22}>0$ and $Q_2=C^T_2C_2\geq0$. It follows that:
\begin{eqnarray*}
  C_2x^*_k=0,  u^{1*}_k=u^{2*}_k=0, N\geq0.
\end{eqnarray*}

The observability of $(A, Q^{\frac{1}{2}}_2)$ given in Assumption 2 indicates that $x_0 = x = 0$, which is a contradiction with $x\neq 0$,
i.e., there exists $N_0$ such that $P^2_0(N) > 0$ for $N>N_0$.

And we will show the uniqueness in the following paragraph.

Since $P^1>0$ satisfies the standard Riccati equation (\ref{opl-eq47-1}), then $P^1$ is uniquely solvable. Suppose $Z^2$ is another solution to (\ref{opl-eq47-1}), (\ref{opl-eq47-2}) and (\ref{opl-eq50-1})-(\ref{opl-eq50-10}) satisfying $Z^2>0$, i.e.,
\begin{eqnarray}
P^1&\hspace{-0.8em}=&\hspace{-0.8em} Q_1+A^TP^1A-A^TP^1B_1(\Gamma^1)^{-1}B^{T}_1P^1A,\label{opl-eq59-1}\\
Z^2&\hspace{-0.8em}=&\hspace{-0.8em}Q_2+A^TY^TR_{21}YA+A^TM^{1T}(\Upsilon^{-1})^TZ^{2}\Upsilon^{-1}M^1A-M^{2T}(\Gamma^2)^{-1} M^2.\label{opl-eq59-2}\\
\Gamma^1&\hspace{-0.8em}=&\hspace{-0.8em}R_{11}+B^T_1P^1B_1,\label{opl-eq59-3}\\
\Gamma^2&\hspace{-0.8em}=&\hspace{-0.8em}R_{22}+B^T_2Y^TR_{21}YB_2+B^T_2M^{1T}(\Upsilon^{-1})^TZ^{2}\Upsilon^{-1}M^1B_2,\label{opl-eq59-4}\\
 M^{1}&\hspace{-0.8em}=&\hspace{-0.8em} I-B_1(\Gamma^1)^{-1}B^T_1P^1,\label{opl-eq59-5}\\
 M^{2}&\hspace{-0.8em}=&\hspace{-0.8em}B^T_2Y^TR_{21}YA+B^T_2M^{1T}(\Upsilon^{-1})^TZ^{2}\Upsilon^{-1}M^1A,\label{opl-eq59-6}\\
 S&\hspace{-0.8em}=&\hspace{-0.8em}(\Gamma^1)^{-1}B^T_1P^1,\label{opl-eq59-7}\\
 \Upsilon&\hspace{-0.8em}=&\hspace{-0.8em}I+B_1(\Gamma^1)^{-1}B^T_1T,\label{opl-eq59-8}\\
 Y &\hspace{-0.8em}=&\hspace{-0.8em}S+(\Gamma^1)^{-1}B^T_1T\Upsilon^{-1}M^1,\label{opl-eq59-9}\\
 K^1&\hspace{-0.8em}=&\hspace{-0.8em}-Y(A+B_2K^2),\label{opl-eq59-10}\\
 K^2&\hspace{-0.8em}=&\hspace{-0.8em}-(\Gamma^2)^{-1} M^2,\label{opl-eq59-11}\\
 T &\hspace{-0.8em}=&\hspace{-0.8em} A^TM^{1T}P^1B_2K^2+A^TM^{1T}T\Upsilon^{-1}M^1(A+B_2K^2).\label{opl-eq59-12}
\end{eqnarray}

According to the  uniqueness of the optimal cost function (\ref{opl-eq51-1}), then we can obtain
\begin{eqnarray}\label{opl-eq60}
  J^*_2 =x^T_0P^2x_0= x^T_0Z^2x_0.
\end{eqnarray}
Since $x_0$ is the arbitrary initial state, thus (\ref{opl-eq60}) implies that $P^2 = Z^2$, the uniqueness of the solution to (\ref{opl-eq47-2}) has been obtained.

\emph{Sufficiency:} Under Assumptions 1 and 2, and if $\Upsilon$ is invertible, suppose $P^1>0$ and $P^2>0$ are the solutions to (\ref{opl-eq47-1}) and (\ref{opl-eq47-2}), respectively, we shall show that $u^1_k\in \mathcal{U}_1$ is to minimize the cost function $J_1$ and $u^2_k\in \mathcal{U}_2$ stabilizes the system (\ref{opl-eq1}) in stabilization and minimizes the cost function $J_2$.

Denote the Lyapunov function $\bar{H}^2_k$ as
\begin{eqnarray}\label{opl-eq61}
  \bar{H}^2_k= x^T_kP^2x_k.
\end{eqnarray}

By using (\ref{opl-eq48}) and (\ref{opl-eq49}), we have
\begin{eqnarray}\label{opl-eq62}
  u^1_k= -YAx_k-YB_2u^2_k,
\end{eqnarray}
adding (\ref{opl-eq62}) into (\ref{opl-eq1}), it derives
\begin{eqnarray}\label{opl-eq63}
x_{k+1}&\hspace{-0.8em}=&\hspace{-0.8em}Ax_k-B_1YAx_k-B_1YB_2u^2_k+B_2u^2_k\nonumber\\
&\hspace{-0.8em}=&\hspace{-0.8em}(I-B_1Y)Ax_k+(I-B_1Y)B_2u^2_k\nonumber\\
&\hspace{-0.8em}=&\hspace{-0.8em}\gamma^{-1}M^1Ax_k+\gamma^{-1}M^1B_2u^2_k.
\end{eqnarray}
Then, combining (\ref{opl-eq61}) and (\ref{opl-eq63}), there follows
\begin{eqnarray}\label{opl-eq64}
&\hspace{-0.8em}&\hspace{-0.8em}\bar{H}^2_k-\bar{H}^2_{k+1}\nonumber\\
=&\hspace{-0.8em}&\hspace{-0.8em}-u^{2T}_k[R_{22}+B^T_2Y^TR_{21}YB_2+B^T_2M^{1T}(\Upsilon^{-1})^TP^{2}\Upsilon^{-1}M^1B_2]u^{2}_k-2u^{2T}_k[B^T_2Y^TR_{21}YA\nonumber\\
&\hspace{-0.8em}&\hspace{-0.8em}+B^T_2M^{1T}(\Upsilon^{-1})^TP^{2}\Upsilon^{-1}M^1A]x_k
+x^T_k[P^2-A^TY^TR_{21}YA-A^TM^{1T}(\Upsilon^{-1})^TP^{2}\Upsilon^{-1}M^1A]x_k\nonumber\\
&\hspace{-0.8em}&\hspace{-0.8em}+u^{1T}_kR_{21}u^1_k+u^{2T}_kR_{22}u^2_k\nonumber\\
=&\hspace{-0.8em}&\hspace{-0.8em}-[u^{2}_k+(\Gamma^2)^{-1}M^2x_k]^T\Gamma^2[u^{2}_k+(\Gamma^2)^{-1}M^2x_k]+x^T_kQ_2x_k+u^{1T}_kR_{21}u^1_k+u^{2T}_kR_{22}u^2_k\nonumber\\
=&\hspace{-0.8em}&\hspace{-0.8em}x^T_kQ_2x_k+u^{1T}_kR_{21}u^1_k+u^{2T}_kR_{22}u^2_k\geq 0,
\end{eqnarray}
where $u^{2}_k=K_2x_k$ for $k\geq 0$ has been imposed in the last identity.

Thus, we can conclude that $\bar{H}^2_k$ is monotonically decreasing with respect to $k$. And Since $P^2>0$, then $\bar{H}^2_k\geq 0$ is bounded below. According to monotone bounded theorem, $\bar{H}^2_k$ is convergent. We will show $\lim_{k\rightarrow\infty} x_k=0$ in the following paragraph.

Actually, let $m$ be any nonnegative integer, by adding from $k = m$ to $k = m+N$ on both sides of (\ref{opl-eq64}) and taking limitation
of $m$, it holds that
\begin{eqnarray}\label{opl-eq65}
 0&\hspace{-0.8em}=&\hspace{-0.8em}\lim_{m\rightarrow\infty}[\bar{H}^2_m-\bar{H}^2_{m+N+1}]\lim_{m\rightarrow\infty}\sum\limits_{k=m}^{m+N}(x^T_kQ_2x_k+u^{1T}_kR_{21}u^1_k+u^{2T}_kR_{22}u^2_k)\geq0.
\end{eqnarray}

Since the coefficient matrices in (\ref{opl-eq1}) and (\ref{opl-eq2-2}) are time invariant, then via a time-shift of length $m$, it yields that
\begin{eqnarray}\label{opl-eq66}
&\hspace{-0.8em}&\hspace{-0.8em}\sum\limits_{k=m}^{m+N}(x^T_kQ_2x_k+u^{1T}_kR_{21}u^1_k+u^{2T}_kR_{22}u^2_k)\geq x^T_mP^2_m(m+N)x_m=x^T_mP^2_0(N)x_m\geq0.
\end{eqnarray}
Taking limitation on both sides of (\ref{opl-eq66}) and using (\ref{opl-eq65}), we know that
\begin{eqnarray}\label{opl-eq67}
 \lim_{m\rightarrow\infty} x^T_mP^2_0(N)x_m=0.
\end{eqnarray}

Moreover, there exists integer $N_0$ such that for any $N >N_0$, $P^2_0(N_0)>0$, thus (\ref{opl-eq67}) implies that
\begin{eqnarray}\label{opl-eq68}
  \lim_{m\rightarrow\infty} x^T_mx_m=0.
\end{eqnarray}
According to (\ref{opl-eq53-1}), we have
\begin{eqnarray}
   &\hspace{-0.8em}&\hspace{-0.8em}\lim_{m\rightarrow\infty}x^T_{m-1}(A+B_1K^1+B_2K^2)^T(A+B_1K^1+B_2K^2)  x_{m-1}=0,
\end{eqnarray}
which implies that
\begin{eqnarray}
\lim_{m\rightarrow\infty}x_{m}=0
\end{eqnarray}
due to the boundedness of $A+B_1K^1+B_2K^2$.

Therefore, the controllers (\ref{opl-eq48}) stabilizes system (\ref{opl-eq1}) in stabilization.

To complete the proof the Theorem 2, we shall show that the optimal controllers (\ref{opl-eq49}) in Theorem 2 minimizes the cost function (\ref{opl-eq2-1}).

Denote the Lyapunov function
\begin{eqnarray}\label{opl-eq69}
 \bar{H}^1_k=x^T_kP^1x_k+x^T_k\zeta_{k-1}.
\end{eqnarray}
Then combining (\ref{opl-eq1}) and (\ref{opl-eq69}), we have
\begin{eqnarray}\label{opl-eq70}
&\hspace{-0.8em}&\hspace{-0.8em}\bar{H}^1_k- \bar{H}^1_{k+1}\nonumber\\
&\hspace{-0.8em}=&\hspace{-0.8em}x^T_kP^1x_k+x^T_k\zeta_{k-1}-x^T_{k+1}P^1x_{k+1}+x^T_{k+1}\zeta_{k}-\zeta^T_kx_{k+1}+\zeta^T_kx_{k+1}\nonumber\\
&\hspace{-0.8em}=&\hspace{-0.8em}u^{1T}_kR_{11}u^1_k+u^{2T}_kR_{12}u^2_k-u^{1T}_k(R_{11}+B^T_1P^1B_1)u^1_k-2u^{1T}_k(B^T_1P^1Ax_k+B^T_1P^1B_2u^2_k+B^T_1\zeta_{k})\nonumber\\
&\hspace{-0.8em}&\hspace{-0.8em}+x^T_k(P^1-A^TP^1A)x_k-2u^{2T}_kB^T_2P^1Ax_k-u^{2T}_k(R_{12}+B^T_2P^1B_2)u^2_k+x^T_k\zeta_{k-1}+\zeta^T_kx_{k+1}\nonumber\\
&\hspace{-0.8em}&\hspace{-0.8em}-x^T_kA^T\zeta_k-u^{2T}_kB^T_2\zeta_k-\zeta^T_kAx_k-\zeta^T_kB_2u^2_k\nonumber\\
&\hspace{-0.8em}=&\hspace{-0.8em}u^{1T}_kR_{11}u^1_k+u^{2T}_kR_{12}u^2_k-[u^1_k+(\Gamma^1)^{-1}(B^T_1P^1Ax_k+B^T_1P^1B_2u^2_k+B^T_1\zeta_{k})]\Gamma^1[u^1_k+(\Gamma^1)^{-1}\nonumber\\
&\hspace{-0.8em}&\hspace{-0.8em}\times (B^T_1P^1Ax_k+B^T_1P^1B_2u^2_k+B^T_1\zeta_{k})]+x^T_k[P^1-A^TP^1A-A^TP^1B_1(\Gamma^1)^{-1}B^{T}_1P^1A]x_k\nonumber\\
&\hspace{-0.8em}&\hspace{-0.8em}-2u^{2T}_k[B^T_2P^1A-B^T_2P^1B_1(\Gamma^1)^{-1}B^T_1P^1A]x_k-u^{2T}_k[R_{12}+B^T_2P^1B_2-B^T_2P^1B_1(\Gamma^1)^{-1}B^T_1P^1\nonumber\\
&\hspace{-0.8em}&\hspace{-0.8em}\times B_2]u^{2}_k+u^{2T}_kB^T_2P^1B_1(\Gamma^1)^{-1}B^T_1\zeta_k+\zeta^T_kB_1(\Gamma^1)^{-1}B^T_1P^1B_2u^2_k+x^T_kA^TP^1B_1(\Gamma^1)^{-1}B^T_1\zeta_k\nonumber\\
&\hspace{-0.8em}&\hspace{-0.8em}+\zeta^T_kB_1(\Gamma^1)^{-1}B^T_1P^1Ax_k+\zeta^T_kB_1(\Gamma^1)^{-1}B^T_1\zeta_k-x^T_kA^T\zeta_k-\zeta^T_kAx_k-u^{2T}_kB^T_2\zeta_k-\zeta^T_kB_2u^2_k\nonumber\\
&\hspace{-0.8em}&\hspace{-0.8em}+x^T_k\zeta_{k-1}+\zeta^T_kx_{k+1}\nonumber\\
&\hspace{-0.8em}=&\hspace{-0.8em}x^T_kQ_1x_k+u^{1T}_kR_{11}u^1_k+u^{2T}_kR_{12}u^2_k-[u^1_k+(\Gamma^1)^{-1}(B^T_1P^1Ax_k+B^T_1P^1B_2u^2_k+B^T_1\zeta_{k})]\Gamma^1\nonumber\\
&\hspace{-0.8em}&\hspace{-0.8em}\times [u^1_k+(\Gamma^1)^{-1}(B^T_1P^1Ax_k+B^T_1P^1B_2u^2_k+B^T_1\zeta_{k})]-2u^{2T}_kB^T_2P^1M^1Ax_k-x^T_kA^TM^{1T}\zeta_k\nonumber\\
&\hspace{-0.8em}&\hspace{-0.8em}-\zeta^T_kM^1Ax_k-u^{2T}_kB^T_2M^{1T}\zeta_k-\zeta^T_kM^1B_2u^2_k-u^{2T}_k(R_{12}+B^T_2P^1M^1B_2)u^{2}_k+x^T_k\zeta_{k-1} \nonumber\\
&\hspace{-0.8em}&\hspace{-0.8em}+\zeta^T_kx_{k+1}+\zeta^T_kB_1(\Gamma^1)^{-1}B^T_1\zeta_k,
\end{eqnarray}
where
\begin{eqnarray}\label{opl-eq71}
-2u^{2T}_kB^T_2P^1M^1Ax_k&\hspace{-0.8em}=&\hspace{-0.8em}-u^{2T}_kB^T_2P^1M^1Ax_k-x^T_kA^TM^{1T}P^1B_2u^2_k\nonumber\\
&\hspace{-0.8em}=&\hspace{-0.8em}-(\zeta_{k-1}-A^TM^{1T}\zeta_k)^Tx_k-x^T_k(\zeta_{k-1}-A^TM^{1T}\zeta_k)\nonumber\\
&\hspace{-0.8em}=&\hspace{-0.8em}-\zeta^T_{k-1}x_k-x^T_k\zeta_{k-1}+x^T_kA^TM^{1T}\zeta_k+\zeta^T_kM^1Ax_k.
\end{eqnarray}

Adding (\ref{opl-eq71}) into (\ref{opl-eq70}), then it yields that
\begin{eqnarray}\label{opl-eq72}
&\hspace{-0.8em}&\hspace{-0.8em}x^T_kQ_1x_k+u^{1T}_kR_{11}u^1_k+u^{2T}_kR_{12}u^2_k\nonumber\\
&\hspace{-0.8em}=&\hspace{-0.8em}\bar{H}^1_k- \bar{H}^1_{k+1}+[u^1_k+(\Gamma^1)^{-1}(B^T_1P^1Ax_k+B^T_1P^1B_2u^2_k+B^T_1\zeta_{k})]\Gamma^1[u^1_k+(\Gamma^1)^{-1}(B^T_1P^1Ax_k\nonumber\\
&\hspace{-0.8em}&\hspace{-0.8em} +B^T_1P^1B_2u^2_k+B^T_1\zeta_{k})]+u^{2T}_k(R_{12}+B^T_2P^1M^1B_2)u^{2}_k+u^{2T}_kB^T_2M^{1T}\zeta_k+\zeta^T_kM^1B_2u^2_k\nonumber\\
&\hspace{-0.8em}&\hspace{-0.8em}-\zeta^T_kB_1(\Gamma^1)^{-1}B^T_1\zeta_k+\zeta^T_{k-1}x_k-\zeta^T_kx_{k+1}.
\end{eqnarray}

Thus, it holds
\begin{eqnarray}\label{opl-eq73}
&\hspace{-0.8em}&\hspace{-0.8em}\sum\limits_{k=0}^N(x^T_kQ_1x_k+u^{1T}_kR_{11}u^1_k+u^{2T}_kR_{12}u^2_k)\nonumber\\
&\hspace{-0.8em}=&\hspace{-0.8em}x^T_0P^1x_0+x^T_0\zeta_{-1}-\bar{H}^1_{N+1}+\sum\limits_{k=0}^N[u^1_k+(\Gamma^1)^{-1}(B^T_1P^1Ax_k+B^T_1P^1B_2u^2_k+B^T_1\zeta_{k})]\Gamma^1\nonumber\\
&\hspace{-0.8em}&\hspace{-0.8em}\times [u^1_k+(\Gamma^1)^{-1}(B^T_1P^1Ax_k +B^T_1P^1B_2u^2_k+B^T_1\zeta_{k})]+\zeta^T_{-1}x_0+\sum\limits_{k=0}^N[u^{2T}_k(R_{12}+B^T_2P^1M^1B_2)u^{2}_k\nonumber\\
&\hspace{-0.8em}&\hspace{-0.8em}+u^{2T}_kB^T_2M^{1T}\zeta_k+\zeta^T_kM^1B_2u^2_k-\zeta^T_kB_1(\Gamma^1)^{-1}B^T_1\zeta_k].
\end{eqnarray}

Noting $P^1>0$ and using (\ref{opl-eq68}), it holds that
\begin{eqnarray}\label{opl-eq74}
\lim_{N\rightarrow\infty}\bar{H}^1_{N+1}=0.
\end{eqnarray}
Taking limitation of $N\rightarrow\infty$ on both sides of (\ref{opl-eq73}) we know the optimal controller minimizes (\ref{opl-eq2-1}) is exactly (\ref{opl-eq49}).

Thus, the optimal cost function of (\ref{opl-eq2-1}) is
\begin{eqnarray}\label{opl-eq76}
  J^*_1= x^T_0(P^1+T+T^T+\Xi)x^T_0,
\end{eqnarray}
which is exactly (\ref{opl-eq51}).

And applying similar procedure to the Lyapunov function of $H^2_k$ in (\ref{opl-eq61}), the controller (\ref{opl-eq48}) minimizes the cost
function (\ref{opl-eq2-2}) can be immediately obtained by (\ref{opl-eq64}), and the optimal cost function is exactly (\ref{opl-eq51-1}).

This completes the proof.
$\hfill\blacksquare$

\section{Numerical Examples}
Consider the system (\ref{opl-eq1}) and the cost functions (\ref{opl-eq2-1}) and (\ref{opl-eq2-2}) with
$A=1$, $B_1=2$, $B_2=1$, $Q_1=1$, $Q_2=1$, $R_{11}=0.6$, $R_{12}=1$, $R_{21}=1$, $R_{22}=6.2$, $x_0=7.3$. By solving (\ref{opl-eq47-1})-(\ref{opl-eq47-2}) and (\ref{opl-eq50-1})-(\ref{opl-eq50-10}), we have $\Upsilon=0.9965>0$, $P^1=1.1325>0$, $P^2=1.2044>0$, $K^1=-0.4268$ and $K^2=-0.0330$. According to Theorem 2, there exist $u^2_k=-0.0330x_k$ stabilizes the system (\ref{opl-eq1}) in stabilization and $u^1_k=-0.4268x_k$ minimize the cost function $J_1$. As shown in Fig. 1, the regulated state is stable.
\begin{figure}[htbp]
\centerline{\includegraphics[width=0.53\textwidth]{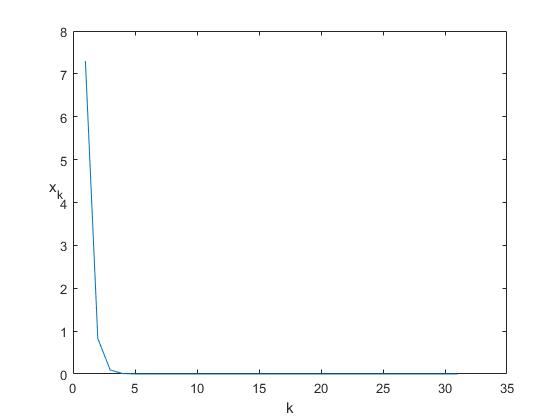}}
\caption{ Trajectory of $x_k$ with controllers $u^1_k=-0.4268 x_k$ and $u^2_k=-0.0330 x_k$.}
\label{fig}
\end{figure}
\section{Conclusion}

In this paper, the stabilization of Stackelberg GBCSs is considered, where the follower is designed to minimise its own cost function and the leader is designed to stabilize the system. To solve the stabilization problem, we consider the finite-time horizon open-loop Stackelberg strategy firstly. Using the maximum principle, we optimize the follower firstly and find that there exist an non-homogeneous relationship between the state and the costate equation. By iteratively solving the FBDEs, the homogeneous relationship between the state and the inhomogeneous terms in costae is obtained. And then by the matrix maximum principle, the optimization of the leader is converted into find the optimal feedback gain matrix, which minimize the cost function of the leader and is subject to the state. Finally, the necessary and sufficient condition for the stabilization generated by the leader is derived.

\end{document}